\documentclass{amsart}%
\usepackage{amsfonts}
\usepackage{amsmath}
\usepackage{amssymb}
\usepackage{graphicx}%
\providecommand{\U}[1]{\protect\rule{.1in}{.1in}}
\newtheorem{theorem}{Theorem}
\theoremstyle{plain}

\newtheorem{corollary}{Corollary}

\newtheorem{definition}{Definition}

\newtheorem{lemma}{Lemma}

\newtheorem{proposition}{Proposition}
\newtheorem{remark}{Remark}

\numberwithin{equation}{section}
\begin{document}
\title[Critical Exponent for Evolution Equations in Modulation Spaces]{Critical Exponent for Evolution Equations in Modulation Spaces}
\author{Huang Qiang}
\address[Huang Qiang]{Department~of~Mathematics,~Zhejiang~University,~Hangzhou~310027,~PR~China}
\email[Huang Qiang]{huangqiang0704@163.com}
\author{Fan~Dashan}
\address[Fan~Dashan]{Department~of~Mathematics,~University~of~Wisconsin-Milwaukee,~Milwaukee,~WI~53201,~USA}
\email[Fan~Dashan]{fan@uwm.edu}
\author{Chen~Jiecheng}
\address[Chen~Jiecheng]{Department~of~Mathematics,~Zhejiang~Normal~University,~Jinhua~321004,~PR~China}
\email{jcchen@zjnu.edu.cn}
\date{November 11, 2014}
\dedicatory{ }
\begin{abstract}
In this paper, we propose a method to find the critical exponent for certain
evolution equations in modulation spaces. We define an index $\sigma (s,q)$,
and use it to determine the critical exponent of the fractional heat
equation as an example. We prove that when $\sigma (s,q)$ is greater than
the critical exponent, this equation is locally well posed in the space $%
C(0,T;M_{p,q}^{s})$; and when $\sigma (s,q)$ is less than the critical
exponent, this equation is ill-posed in the space $C(0,T;M_{2,q}^{s})$. Our
method may further be applied to some other evolution equations.

\end{abstract}
\maketitle

\section{Introduction and main results}

\hspace{6mm} As we all know, many evolution equations have their critical
exponents on either Sobolev spaces or Besov spaces, or both. For example,
the critical exponent of nonlinear Schr\"{o}dinger equation (NLS) in Besov
spaces $\dot{B}_{p,2}^{s}$ is $\frac{n}{p}-\frac{2}{k-1}$ where $k$ is the
power of the nonlinear term \ $u^{k}.$ Cazenave and Weissler \cite{CZ}
showed that NLS is locally well-posed in $C([-T,T];\dot{H}^{s})$ when $s\geq
0$ and $s\geq \frac{n}{2}-\frac{2}{k-1}$. In \cite{CCT}, Christ, Colliander
and Tao proved that when $s<\max \{0,\frac{n}{2}-\frac{2}{k-1}\}$, NLS is
ill-posed in $\dot{H}^{s}$. In \cite{MXZ}, Miao, Xu and Zhao proved similar
results for the nonlinear Hartree equation. We observe that both works in
\cite{CCT} and \cite{MXZ} are heavily based on the scaling invariance of the
work space. On the other hand, the modulation space \ $M_{p,q}^{s}$ \ is
lack of the scaling property, although this space emerges in recent years
and plays a significant role in the study of certain nonlinear evolution
equations. (We will describe more details of the modulation space in the
following content.) \ Since we are not able to find in literature any study
on critical exponent for evolution equation in the modulation space, the aim
of this paper is to propose a different method from \cite{CCT} and \cite{MXZ}%
\ to find the critical exponents. Particularly we find the critical exponent
for the fractional heat equation on the modulation space, without the
scaling invariance. This exponent satisfies the well and ill posedness
property on the modulation space, which is quite similar to that for NLS in
the Sobolev space.

Modulation spaces was introduced by Feichtinger in \cite{F} to measure
smoothness of a function or distribution in a way different from $L^{p}$
spaces, and they are now recognized as a useful tool for studying
pseudo-differential operators (see \cite{B}\cite{CFS}\cite{M}\cite{S}\cite%
{Tj}). The original definition of the modulation space is based on the
short-time Fourier transform and window function. In \cite{WHH}, Wang and
Hudizk gave an equivalent definition of the discrete version on modulation
spaces by the frequency-uniform-decomposition. With this discrete version,
they are able to study the global solution for nonlinear Schr\"{o}dinger
equation and nonlinear Klein-Gordon equation. After then, there are many
studies on nonlinear PDEs in modulation spaces followed their work. Below we
list some of them, among many others. In \cite{GC}, Guo and Chen proved the
Stricharz estimates on $\alpha $-modulation spaces. For well-posedness in
modulation space, Wang, Zhao, and Guo \cite{WZG} studied the local solution
for nonlinear Schr\"{o}dinger equation and Navier-Stokes equations. In \cite%
{WH}, Wang and Huang studied the local and global solutions for generalized
KdV equations, Benjamin-Ono and Schr\"{o}dinger equations. In \cite{RMW},
Ruzhansky, Sugimoto and Wang stated some new progress and open questions in
modulation spaces. Also, for the ill posedness in modulation spaces,
Iwabuchi studied well and ill posedness for Navier-Stokes equations and heat
equations (see \cite{TI}). Iwabuchi's result can be stated in the following
theorem:

\textbf{Theorem A} \bigskip \cite{TI}:When $s-\frac{n}{q^{\prime }}>-\frac{2}{k-1}$,
the Heat equation
\begin{equation*}
(H)~~~~u(t)=e^{t\Delta }u_{0}+\int_{0}^{t}e^{t-\tau \Delta }u^{k}d\tau
\end{equation*}%
is locally well-posed in $C([0,T,M_{p,q}^{s}])$. When $s<-\frac{2}{k}$ or $s-%
\frac{n}{q^{\prime }}<-\frac{n+2}{k}$, equation (H) is ill-posed in $%
C([0,T,M_{2,q}^{s}])$

Since Iwabuchi's result is not a sharp one, a natural question is if there
are some critical exponents for this equation in modulation spaces based on
the well and ill posedness. In this paper, we will answer this question.

\bigskip

First, we recall some important properties of Besov spaces \cite{T}. The
first one is a Sobolev embedding that says $B_{p_{1},q}^{s_{1}}\subset
B_{p_{2},q}^{s_{2}}$ if and only if
\begin{equation*}
s_{2}\leq s_{1}~~and~~s_{1}-\frac{n}{p_{1}}=s_{2}-\frac{n}{p_{2}}.
\end{equation*}%
The second one says that the Besov space $B_{p,q}^{s}$ \ forms a
multiplication algebra if $s-\frac{n}{p}>0$. By comparing these properties
to the algebra property of modulation spaces and (2.2) in Section 2, we
observe that the index $s-\frac{n}{p}$ in the Besov space is an analog of
the index $s-n(1-\frac{1}{q})$ in the modulation space. Motivated by such an
observation, heuristically, we may use the index $\ s-n(1-\frac{1}{q})$ to
describe the critical exponent in the modulation space. Of course, this
heuristic idea will be technically supported in our following discussion.
For convenience in the discussion, we denote $\sigma (s,q)=s-n(1-\frac{1}{q}%
) $, and use the inequality
\begin{equation*}
A(u,v,w...)\preceq B(u,v,w...)
\end{equation*}%
to mean that there is a positive number $C$ independent of all main
variables $\ u,v,w...$, for which $A(u,v,w...)\leq CB(u,v,w...)$.

\bigskip

Now we state a general theorem for well posedness.

\begin{theorem}
Let $U(t)$ be the dispersive semigroup:
\begin{equation*}
U(t):=\mathcal{F}^{-1}e^{tp(\xi )}\mathcal{F}
\end{equation*}%
where $p(\xi ):R^{n}\rightarrow C$, $\mathcal{F}$ denotes the Fourier
transform. Assume that there exists a $\ \theta >0$ \ such that $U(t)$
satisfies the estimate: \ for \ $1\leq p<\infty ,$
\begin{equation*}
\Vert U(t)f\Vert _{M_{p,q_{1}}^{s_{1}}}\preceq t^{-\theta \lbrack \sigma
(s_{1},q_{1})-\sigma (s_{2},q_{2})]}\Vert f\Vert _{M_{p,q_{2}}^{s_{2}}}\eqno%
(1.1)
\end{equation*}%
for all $\ -\infty <s_{2}\leq s_{1}<+\infty ,$ $\ 1\leq q_{1}\leq
q_{2}<+\infty $ \ and \ $0<t<\infty $. Then the general dispersive equation
\begin{equation*}
u=U(t)u_{0}+\int_{0}^{t}U(t-\tau )u^{k}d\tau ,~~~~~~~~k\in Z^{+}\eqno(1.2)
\end{equation*}%
is locally and globally well-posed for any $s>0,q\geq 1$ and $\sigma (s,q)>-%
\frac{1}{(k-1)\theta }$. More precisely, we have the following statements.%
\newline
(i) \ Let $1\leq p,q<\infty ,s>0$ and $\sigma (s,q)>-\frac{1}{(k-1)\theta }$%
. For any $u_{0}\in M_{p,q}^{s}$, there exists a $\ T>0$ such that the
equation (1.2) has an unique solution in $C(0,T;M_{p,q}^{s})$.\newline
(ii) There exists a small number $\nu >0$ such that for any $\Vert
u_{0}\Vert _{M_{2,2}^{s}}\leq \nu $, the equation (1.2) has an unique
solution in the space
\begin{equation*}
L^{\infty }(R;M_{2,2}^{s})\bigcap\limits_{1\leq q<2}L^{\gamma
(q)}(R;M_{2,q}^{s}),
\end{equation*}%
where $\frac{1}{\gamma (q)}=(\frac{1}{q}-\frac{1}{2})n\theta <\frac{1}{2}$, $%
s>0$ and $\sigma (s,2)>-\frac{1}{(k-1)\theta }.$
\end{theorem}

\begin{remark}
In this theorem, we can see that the index $q^{\prime }$ in modulation
spaces plays a similar role as the index $p$ in Besov spaces (see \cite{T}).
The only difference is that $s-\frac{n}{2}$ can not be equal to $-\frac{1}{%
(k-1)\theta }$ in the global case. This is because that the equality does
not hold in the condition in  (2.2), since Wang and Hudzik in \cite{WHH}
proved that the condition in (2.2) is sharp. The reader can find this
condition in Section 2.
\end{remark}

\begin{remark}
It is well known that the Schr\"{o}dinger semigroup $S(t)$ has the following
estimate in Besov spaces for $2\leq p<\infty $:
\begin{equation*}
\Vert S(t)f\Vert _{B_{p,2}^{s}}\preceq t^{-(\frac{1}{2}-\frac{1}{p})n}\Vert
f\Vert _{B_{p^{\prime },2}^{s}}\eqno(1.3)
\end{equation*}%
If we rewrite above inequality as following:
\begin{equation*}
\Vert S(t)f\Vert _{B_{p,2}^{s}}\preceq t^{-\frac{1}{2}[(s-\frac{n}{p})-(s-%
\frac{n}{p^{\prime }})]}\Vert f\Vert _{B_{p^{\prime },2}^{s}}\eqno(1.4)
\end{equation*}%
we find that its critical exponent is $s_{c}=\frac{n}{p}-\frac{2}{k-1}$.
From this observation we see that $\ \theta $ \ in Theorem 1 just likes $%
\frac{1}{2}$ in (1.4).
\end{remark}

Now, as an application of Theorem 1, we consider the Cauchy problem for
the fractional heat equation
\begin{equation*}
u_{t}+\left( -\Delta \right) ^{\frac{\alpha }{2}}u=u^{k},~~~~~~~u(0)=u_{0}%
\eqno(1.5)
\end{equation*}%
The following two theorems show that in this equation, $\theta =\frac{1}{%
\alpha }$ is the minimum number in the inequality (1.1) for the fractional
heat equation,\ and \ $-\frac{\alpha }{k-1}$ \ is critical.

\begin{theorem}
Let $1\leq p,q<\infty $, $s\geq 0$ and $\sigma (s,q)>-\frac{\alpha }{k-1}$.
There exists a $\ T>0$ \ such that the equation (1.5) is locally well-posed
in $C(0,T;M_{p,q}^{s}).$
\end{theorem}

\begin{theorem}
Let $1\leq q<\infty $. When $\sigma (s,q)<-\frac{\alpha }{k-1}$ or $s<-\frac{%
\alpha }{k-1}$ for any $q$, then there exists a $\ T>0$ \ for which the
equation (1.5) is ill-posed in $C(0,T;M_{2,q}^{s}).$
\end{theorem}




Comparing above results (Theorem 2, Theorem 3) in the case \ $\alpha =1$ \ and Iwabuchi's result (Theorem A). For the area $s>0$, our result is a sharp one
which Theorem A is not. On the other hand, our result not only works for \ $\alpha =1,$ but it gives the critical exponent of (1.5)
for all \ $\alpha >0.$ By the same way, We can get the similar sharp result for incompressible Navier-Stokes equations which is also 
better than Iwabuchi's (see \cite{TI})

It is known that for each evolution equation, we have a set \ $\Theta $ \ of
indices \ $\theta $ \ for which the time-spaces estimate (1.1) holds for \ $%
0<t\leq 1$, and we have a critical exponent for its Cauchy problem with
nonlinear term \ $u^{k}.$ It is reasonable to guess that if the set \ $%
\Theta $ \ has the positive minimum value and if we obtain the critical
exponent
\begin{equation*}
\sigma (s,q)>-\frac{1}{(k-1)\theta _{0}}
\end{equation*}%
on the modulation space \ $M_{p,q}^{s},$ then this \ $\theta _{0}$ \ must be
the minimum value of \ $\Theta .$ Although these results seem to work for
the case $\ \theta _{0}>0,$ if  $p(\xi )$ in the symbol of the fundamental
semi-group is real, the method used in our proof may also work in the case \
$\theta _{0}=0.$ For instance, we look the Schr\"{o}dinger equation, for
which we can not obtain the time-space estimate as (1.1). Actually we only
have the estimate:
\begin{equation*}
\Vert S(t)f\Vert _{M_{2,q}^{s}}\leq \Vert f\Vert _{M_{2,q}^{s}}.\eqno(1.6)
\end{equation*}%
We can still use the same method to obtain partial conclusion as that for
the fractional heat equation. The following theorem is our result for Schr%
\"{o}dinger equation:

\begin{corollary}
Let $1\leq q<\infty $. When $s>0$ and $\sigma (s,q)>0$, the Schr\"{o}dinger
equation is locally well-posed in $C(0,T;M_{2,q}^{s})$ for some \ $T>0$.
When $s=0$ and $\sigma (s,q)\geq 0$, the Schr\"{o}dinger equation is locally
well-posed in $C(0,T;M_{2,q}^{0})$ for some \ $T>0$. When $\sigma (s,q)<-%
\frac{2}{k-1}$ or $s<-\frac{2}{k-1}$, the Schr\"{o}dinger equation is
ill-posed in $C(0,T;M_{2,q}^{s})$.
\end{corollary}

Also, for Klein-Gordon equation, if we write in this form:
\begin{equation*}
u(t)=K(t)(u_{0}+\frac{tant\omega^{\frac{1}{2}}}{\omega^{\frac{1}{2}}}%
u_{1})+\int_{0}^{t}\frac{K_{1}(t-\tau)}{\omega^{\frac{1}{2}}}u^{k}d\tau\eqno%
(1.7)
\end{equation*}
where
\begin{equation*}
K(t)=cost\omega^{\frac{1}{2}},~~~K_{1}(t)=sint\omega^{\frac{1}{2}%
},~~~\omega=(I-\Delta),
\end{equation*}
then choose $u_{0}=\frac{tant\omega^{\frac{1}{2}}}{\omega^{\frac{1}{2}}}u_{1}
$ in the proof of ill posedness, we can obtain following corollary:

\begin{corollary}
Let $1\leq q<\infty $. When $s\geq 0$ and $\sigma (s,q)>-\frac{1}{k-1}$, the
equation (1.7) is locally well-posed in $C(0,T;M_{2,q}^{s})$. When $\sigma
(s,q)<-\frac{2}{k-1}$ or $s<-\frac{2}{k-1}$ for any $q$, then the equation
(1.5) is ill-posed in $C(0,T;M_{2,q}^{s}).$
\end{corollary}

In \ Corollary 1, there is a gap in the interval $\ \ [-\frac{2}{k-1},0]$ \
for the Schr\"{o}dinger equation. This is an unsolved problem. Similar gaps
exsit for the Klein-Gordon equation in Corollary 2.\

It is interesting to see that the index  $q'$  plays a crucial rule in the
study of modulation space \ $M_{p,q}^{s},$ while it plays almost no role in
the study of the Besov space \ $B_{p,q}^{s}.$ The essence of this phenomenon
is that they have different geometric regions in decompositions on the
frequency space, so that the Bernstein inequality gives quite different
estimates in the proofs of their embedding and algebra properties.

This paper organized as follows. In Section 2, we will introduce some basic
knowledge on the modulation space, as well as some useful estimates that
will be used in our proofs. All proofs of main theorems will be presented in
Section 3.

\section{Preliminaries}

\hspace{6mm}In this section, we give the definition and discuss some basic
properties of modulation spaces. Also, we will prove some estimates which
are descried by the index $\sigma (s,q)$.

\begin{definition}
(Modulation spaces) Let $\{\varphi _{k}\}\subset C_{0}^{\infty }(R^{n})$ be
a partition of the unity satisfying the following conditions:
\begin{equation*}
supp\varphi \subset \{\xi \in R^{n}:\mid \xi \mid \leq \sqrt{n}\},\sum_{k\in
Z^{n}}\varphi (\xi -k)=1,\varphi _{k}(\xi ):=\varphi (\xi -k)
\end{equation*}%
for any $\xi \in R^{n}$, and let
\begin{equation*}
\Box _{k}:=\mathcal{F}^{-1}\varphi _{k}\mathcal{F}.
\end{equation*}
\end{definition}

By this frequency-uniform decomposition operator, we define the modulation
spaces \ $M_{p,q}^{s}(R^{n}),$ for $0<p,q\leq \infty ,$ $-\infty <s<\infty ,$
by
\begin{equation*}
M_{p,q}^{s}(R^{n}):=\{f\in S^{\prime }:\Vert f\Vert
_{M_{p,q}^{s}(R^{n})}=(\sum_{k\in Z^{n}}<k>^{sq}\Vert \Box _{k}f\Vert
_{p}^{q})^{\frac{1}{q}}<\infty \},
\end{equation*}%
where $\langle k\rangle=\sqrt{1+|k|^{2}}$. See \cite{WHH} for details.

\begin{proposition}
(Isomorphism)\cite{WHH} Let $0<p,q\leq \infty ,s,\sigma \in R$. $J_{\sigma
}=(I-\triangle )^{\frac{\sigma }{2}}:M_{p,q}^{s}\rightarrow
M_{p,q}^{s-\sigma }$ is an isomorphic mapping, where \ $I$ \ is the identity
mapping and \ $\Delta $ \ is the Laplacian.
\end{proposition}

\begin{proposition}
(Embedding).\cite{WHH} We have\newline
(i) $M^{s_{1}}_{p_{1},q_{1}}\subset M^{s_{2}}_{p_{2},q_{2}}$,~if $s_{1}\geq
s_{2},0<p_{1}\leq p_{2},0<q_{1}\leq q_{2}.$\qquad\qquad\qquad\qquad\qquad(2.1)\\
(ii) $M^{s_{1}}_{p_{1},q_{1}}\subset M^{s_{2}}_{p_{2},q_{2}}$,~if $%
q_{1}>q_{2},s_{1}>s_{2},s_{1}-s_{2}>n/q_{2}-n/q_{1}.$\qquad\qquad\qquad\qquad(2.2)
\end{proposition}

\begin{lemma}
Let $s\geq0$, $k\in Z^{+}$, $1\leq p,q,q_{1},q_{2}\leq\infty$, and $\frac{1}{%
q}+k-1=\frac{1}{q_{1}}+\frac{k-1}{q_{2}}$. We have
\begin{equation*}
\|u^{k}\|_{M_{p,q}^{s}}\preceq\|u\|_{M_{p,q_{1}}^{s}}\|u%
\|_{M_{p,q_{2}}}^{k-1}.\eqno(2.3)
\end{equation*}
\end{lemma}

\textbf{Proof: } We only consider the case when $k=2$ for simplicity, since
the proof for \ $k\geq 2$ \ is similar. Note the project operators $\left\{
\square _{k}\right\} $ \ satisfying
\begin{equation*}
\sum_{k\in
\mathbb{Z}
^{n}}\square _{k}=I.
\end{equation*}%
By the Minkowski inequality, we may write
\begin{equation*}
\langle i\rangle^{s}\Vert \Box _{i}u^{2}\Vert _{L^{p}}\leq
\langle i\rangle^{s}\sum\limits_{i_{1},i_{2}\in Z^{n}}\Vert \Box _{i}(\Box _{i_{1}}u\Box
_{i_{2}}u)\Vert _{L^{p}}.
\end{equation*}

We observe that the support condition of \ $\square _{k}$ in the frequency
space$\mathcal{\ }$implies that
\begin{equation*}
\Box _{i}(\Box _{i_{1}}u\Box _{i_{2}}u)=0~~if~~|i-i_{1}-i_{2}|\geq k_{0}(n),
\end{equation*}%
where $k_{0}(n)$ is an integer which depends only on $n$ (see \cite{WHH}).
So we have
\begin{equation*}
\langle i\rangle^{s}\Vert \Box _{i}u^{2}\Vert _{L^{p}}\leq
\langle i\rangle^{s}\sum\limits_{i_{1},i_{2}\in Z^{n},|i-i_{1}-i_{2}|\leq k_{0}(n)}\Vert
\Box _{i}(\Box _{i_{1}}u\Box _{i_{2}}u)\Vert _{L^{p}}.
\end{equation*}%
By the Bernstein and H\"{o}lder's inequalities, we obtain that
\begin{eqnarray*}
\langle i\rangle^{s}\Vert \Box _{i}u^{2}\Vert _{L^{p}} &\preceq
&\sum\limits_{i_{1},i_{2}\in Z^{n},|i-i_{1}-i_{2}|\leq
k_{0}(n)}\langle i_{1}+i_{2}\rangle^{s}\Vert \Box _{i_{1}}u\Vert _{L^{p_{1}}}\Vert \Box
_{i_{2}}u\Vert _{L^{p_{2}}} \\
&\preceq &\sum\limits_{i_{1},i_{2}\in Z^{n},|i-i_{1}-i_{2}|\leq
k_{0}(n)}\langle i_{1}\rangle^{s}\Vert \Box _{i_{1}}u\Vert _{L^{p_{1}}}\Vert \Box
_{i_{2}}u\Vert _{L^{p_{2}}} \\
&~&~~~~+\sum\limits_{i_{1},i_{2}\in Z^{n},|i-i_{1}-i_{2}|\leq
k_{0}(n)}\langle i_{2}\rangle^{s}\Vert \Box _{i_{1}}u\Vert _{L^{p_{1}}}\Vert \Box
_{i_{2}}u\Vert _{L^{p_{2}}},
\end{eqnarray*}%
where $\frac{1}{p}=\frac{1}{p_{1}}+\frac{1}{p_{2}}$. Thus, by (2.1) and
Young's inequality of series, we have
\begin{equation*}
\Vert u^{2}\Vert _{M_{p,q}^{s}}\preceq \Vert u\Vert _{M_{p,q_{1}}^{s}}\Vert
u\Vert _{M_{p,q_{2}}}.\eqno(2.4)
\end{equation*}%
By the induction and (2.4), we can easily obtain the desired result.

\begin{remark}
In \cite{C}, Cazenave proved $\Vert u^{k}\Vert _{B_{p,2}^{s}}\preceq \Vert
u\Vert _{B_{p_{1},2}^{s}}\Vert u\Vert _{L^{p_{2}}}^{k-1}$ when $\frac{1}{p}=%
\frac{1}{p_{1}}+\frac{k-1}{p_{2}}$. Also, we can see that the condition of
Lemma 1 is equivalent to $\frac{1}{q^{\prime }}=\frac{1}{q_{1}^{\prime }}+%
\frac{k-1}{q_{2}^{\prime }}.$ This again indicates that the index $\frac{1}{%
q^{\prime }}$ in modulation spaces behaves like the index $\frac{1}{p}$ in
the Besov spaces.
\end{remark}

\begin{lemma}
Let $1\leq p\leq\infty$, $1\leq q_{1}\leq q_{2}$, $0\leq s_{2}\leq s_{1}$
and $\sigma(s_{1},q_{1})-\sigma(s_{2},q_{2})=R.$ If $s_{1}>0$ and $%
\sigma(s_{1},q_{1})>-\frac{R}{k-1}$, then we have
\begin{equation*}
\|u^{k}\|_{M_{p,q_{2}}^{s_{2}}}\preceq\|u\|^{k}_{M_{p,q_{1}}^{s_{1}}}.\eqno%
(2.5)
\end{equation*}
\end{lemma}

\textbf{Proof: } Fix a small \ $\epsilon >0$ \ and pick \ $q_{3}\ \ $such
that
\begin{equation*}
\frac{n}{q_{3}}=\frac{n}{q_{1}}-R+\epsilon .
\end{equation*}%
It is easy to check that $s_{1}+\frac{n}{q_{3}}>s_{2}+\frac{n}{q_{2}}$.
Using (2.2), we have
\begin{equation*}
\Vert u^{k}\Vert _{M_{p,q_{2}}^{s_{2}}}\preceq \Vert u^{k}\Vert
_{M_{p,q_{3}}^{s_{1}}}.\eqno(2.6)
\end{equation*}%
Since $s_{1}>0$, by Lemma (2.3), we obtain that
\begin{equation*}
\Vert u^{k}\Vert _{M_{p,q_{3}}^{s_{1}}}\preceq \Vert u\Vert
_{M_{p,q_{1}}^{s_{1}}}\Vert u\Vert _{M_{p,q_{5}}}^{k-1},\eqno(2.7)
\end{equation*}%
where
\begin{equation*}
\frac{n}{q_{5}}=n-\frac{R}{k-1}+\frac{\epsilon }{k-1},\eqno(2.8)
\end{equation*}%
and $\epsilon $ is small enough to ensure $\ s_{1}+\frac{n}{q_{1}}>\frac{n}{%
q_{5}}$. Using (2.2) again, we have
\begin{equation*}
\Vert u\Vert _{M_{p,q_{5}}}\preceq \Vert u\Vert _{M_{p,q_{1}}^{s_{1}}}.\eqno%
(2.9).
\end{equation*}%
Inserting (2.9) into (2.7), we now obtain (2.5). This completes the proof.

\section{Proof of the main theorems}

\textbf{Proof of Theorem 1 }We first prove the local case. Consider the
integral equation
\begin{equation*}
\Phi (u)=U(t)u_{0}+\int_{0}^{t}U(t-\tau )u^{k}d\tau .
\end{equation*}%
It is well known that this equation is equivalent to the Cauchy problem
(1.5). To prove the above equation has a unique solution, we will use the
standard contraction method. To this end, we define the space
\begin{equation*}
X_{1}=\{u:\Vert u\Vert _{L^{\infty }(0,T;M_{p,q}^{s})}\leq C_{0}\}
\end{equation*}%
with the metric \
\begin{equation*}
~~~~d(u,v)=\Vert u-v\Vert _{L^{\infty }(0,T;M_{p,q}^{s})},
\end{equation*}%
where the positive numbers $C_{0}$ and $T$ will be chosen later when we
invoke the contraction. We now choose numbers $\overline{s}$ $\ $and $\
\overline{q}$ \ for which
\begin{equation*}
\sigma (s,q)-\sigma (\overline{s},\overline{q})=\frac{1}{\theta +\varepsilon
},
\end{equation*}%
where $\varepsilon \ $is a small positive number such that $\sigma (s,q)>-%
\frac{1}{(\theta +\varepsilon )(k-1)}$. By (1.1) and Lemma 2, we have
\begin{eqnarray*}
~~~~~~~~\Vert \Phi (u)\Vert _{X_{1}} &\preceq &\Vert u_{0}\Vert
_{M_{p,q}^{s}}+\Vert \int_{0}^{t}U(t-\tau )u^{k}d\tau \Vert _{X_{1}} \\
&\preceq &\Vert u_{0}\Vert _{M_{p,q}^{s}}+\sup_{t\in (0,T]}\left\vert
\int_{0}^{t}(t-\tau )^{-\theta \lbrack \sigma (s,q)-\sigma (\overline{s},%
\overline{q})]}\Vert u^{k}\Vert _{M_{p,\overline{q}}^{\overline{s}}}d\tau
\right\vert \\
&\preceq &\Vert u_{0}\Vert _{M_{p,q}^{s}}+\sup_{t\in (0,T]}\left\vert
\int_{0}^{t}(t-\tau )^{-\theta \lbrack \sigma (s,q)-\sigma (\overline{s},%
\overline{q})]}d\tau \right\vert \Vert u\Vert _{X_{1}}^{k} \\
&\preceq &\Vert u_{0}\Vert _{M_{p,q}^{s}}+T^{1-\frac{\theta }{\theta
+\varepsilon }}\Vert u\Vert
_{X_{1}}^{k}.~~~~~~~~~~~~~~~~~~~~~~~~~~~~~~~~~~~~~~~~~~~~(3.2)
\end{eqnarray*}%
By the contraction mapping argument, we obtain (i) in Theorem 1 after
choosing suitable \ $T$ \ and \ $C_{0}$.

Next, we consider the global case. Choosing $s_{1}=s_{2},~\frac{1}{q_{1}}+%
\frac{1}{q_{2}}=1$ in (1.1), we have
\begin{equation*}
\Vert U(t)f\Vert _{M_{2,q}^{s_{{}}}}\preceq t^{-\theta \lbrack n(\frac{2}{q}%
-1)]}\Vert f\Vert _{M_{2,q^{\prime }}^{s}}\eqno(3.3)
\end{equation*}%
for $1\leq q\leq 2$, $s>0$. When $\theta \lbrack n(\frac{2}{q}-1)]<1$, we
can obtain the following estimates by standard dual methods (see \cite{WHCHC}%
):
\begin{equation*}
\Vert U(t)f\Vert _{L^{\frac{2}{\theta \lbrack n(\frac{2}{q}-1)]}%
}(R;M_{2,q}^{s})}\preceq \Vert f\Vert _{M_{2,2}^{s}},\eqno(3.4)
\end{equation*}%
\begin{equation*}
\Vert \int_{0}^{t}U(t-\tau )fd\tau \Vert _{L^{\frac{2}{\theta \lbrack n(%
\frac{2}{q}-1)]}}(R;M_{2,q}^{s})}\preceq \Vert f\Vert
_{L^{1}(R;M_{2,2}^{s})},\eqno(3.5)
\end{equation*}%
\begin{equation*}
\Vert \int_{0}^{t}U(t-\tau )fd\tau \Vert _{L^{\infty
}(R;M_{2,2}^{s})}\preceq \Vert f\Vert _{L^{(\frac{2}{\theta \lbrack n(\frac{2%
}{q}-1)]})^{\prime }}(R;M_{2,q^{\prime }}^{s})},\eqno(3.6)
\end{equation*}%
\begin{equation*}
\Vert \int_{0}^{t}U(t-\tau )fd\tau \Vert _{L^{\frac{2}{\theta \lbrack n(%
\frac{2}{q}-1)]}}(R;M_{2,q}^{s})}\preceq \Vert f\Vert _{L^{(\frac{2}{\theta
\lbrack n(\frac{2}{q}-1)]})^{\prime }}(R;M_{2,q^{\prime }}^{s})}.\eqno(3.7)
\end{equation*}%
By interpolation among (3.5), (3.6) and (3.7), we obtain that
\begin{equation*}
\Vert \int_{0}^{t}U(t-\tau )fd\tau \Vert _{L^{\frac{2}{\theta \lbrack n(%
\frac{2}{q}-1)]}}(R;M_{2,q}^{s})}\preceq \Vert f\Vert _{L^{(\frac{2}{\theta
\lbrack n(\frac{2}{r}-1)]})^{\prime }}(R;M_{2,r^{\prime }}^{s})}\eqno(3.8)
\end{equation*}%
for any $1\leq q,r\leq 2$.

We choosing $\frac{1}{q}=\frac{1}{2}+\frac{1}{(k-1)n\theta }$ and let
\begin{equation*}
X_{2}=L^{\infty }(R;M_{2,2}^{s})\bigcap L^{\frac{2}{\theta \lbrack n(\frac{2%
}{q}-1)]}}(R;M_{2,q}^{s})
\end{equation*}%
with the metric
\begin{equation*}
d(u,v)=\Vert u-v\Vert _{L^{\infty }(R;M_{2,2}^{s})}+\Vert u-v\Vert _{L^{%
\frac{2}{\theta \lbrack n(\frac{2}{q}-1)]}}(R;M_{2,q}^{s})}
\end{equation*}%
by (1.1),(3.4),(3.8) and Lemma 2, we have
\begin{eqnarray*}
~~~~~~~~~~~~~~~\Vert \Phi (u)\Vert _{X_{2}} &\preceq &\Vert u_{0}\Vert
_{M_{2,2}^{s}}+\Vert u^{k}\Vert _{L^{(\frac{2}{\theta \lbrack n(\frac{2}{q}%
-1)]})^{\prime }}(R;M_{2,q^{\prime }}^{s})} \\
&\preceq &\Vert u_{0}\Vert _{M_{2,2}^{s}}+\Vert u\Vert _{L^{k(\frac{2}{%
\theta \lbrack n(\frac{2}{q}-1)]})^{\prime }}(R;M_{2,q}^{s})}^{k} \\
&=&\Vert u_{0}\Vert _{M_{2,2}^{s}}+\Vert u\Vert
_{X_{2}}^{k}.\qquad\qquad\qquad\qquad\qquad\qquad\qquad\qquad\qquad(3.9)
\end{eqnarray*}%
it is easy to check that
\begin{equation*}
k(\frac{2}{\theta \lbrack n(\frac{2}{q}-1)]})^{\prime }=\frac{2}{\theta
\lbrack n(\frac{2}{q}-1)]}\eqno(3.10)
\end{equation*}%
and
\begin{equation*}
s+\frac{n}{q}>n-\frac{\frac{n}{q}-\frac{n}{q^{\prime }}}{k-1}.\eqno(3.11)
\end{equation*}%
From (3.10), we have
\begin{equation*}
\frac{nk}{q}-\frac{n}{q^{\prime }}=\frac{1}{\theta }-\frac{n}{2}+\frac{nk}{2}%
,\eqno(3.12)
\end{equation*}%
then insert (3.12) into (3.11), we can obtain
\begin{equation*}
s>\frac{n}{2}-\frac{1}{(k-1)\theta }.\eqno(3.13)
\end{equation*}%
Using the standard contraction mapping argument in (3.9), we can find unique
solution in $X_{2}$. Then by (3.4) and (3.8), we can obtain the conclusion
of \ (ii) in Theorem 1. \newline
\textbf{Proof of Theorem 2. }We first prove
\begin{equation*}
\Vert e^{-t(-\Delta )^{\frac{\alpha }{2}}}f\Vert _{M_{p,q}^{s_{1}}}\preceq
(1+t^{-\frac{1}{\alpha }(s_{1}-s_{2})})\Vert f\Vert _{M_{p,q}^{s_{2}}}\eqno%
(3.14)
\end{equation*}%
for any $s_{1}\geq s_{2}$. For the low frequency part $|k|\leq 100\sqrt{n}$,
we have
\begin{eqnarray*}
\sum\limits_{|k|\leq 100\sqrt{n}} &~&\langle k\rangle^{s_{1}q}\Vert \Box _{k}e^{-t(-\Delta
)^{\frac{\alpha }{2}}}f\Vert _{L^{p}}^{q} \\
&\preceq &\sum\limits_{|k|\leq 100\sqrt{n}}\langle k\rangle^{s_{2}q}\Vert \Box _{k}f\Vert
_{L^{p}}^{q}\preceq \Vert f\Vert _{M_{p,q}^{s_{2}}}^{q}.
\end{eqnarray*}%
For the high frequency part, note that the operator \ $\Box
_{k}e^{-t(-\Delta )^{\frac{\alpha }{2}}}$ \ can be written as
\begin{equation*}
\Box _{k}e^{-t(-\Delta )^{\frac{\alpha }{2}}}=\sum_{\left\vert \ell
\right\vert \leq 1}\Box _{k+\ell }e^{-t(-\Delta )^{\frac{\alpha }{2}}}\Box
_{k}
\end{equation*}%
and \ $\Box _{k+\ell }e^{-t(-\Delta )^{\frac{\alpha }{2}}}$ are convolution
operators with the kernels
\begin{equation*}
\Omega _{k+\ell }(y)=e^{i<k+\ell ,y>}\int_{%
\mathbb{R}
^{n}}e^{-t\left\vert \xi +k+\ell \right\vert ^{\alpha }}e^{i<y,\xi >}\varphi
(\xi )d\xi .
\end{equation*}%
Hence, when $|k|\geq 100\sqrt{n}$ \ it is easy to prove
\begin{equation*}
\Vert \Box _{k}e^{-t(-\Delta )^{\frac{\alpha }{2}}}f\Vert _{L^{p}}\preceq
e^{-\frac{t}{2}|k|^{\alpha }}\Vert \Box _{k}f\Vert _{L^{p}}.
\end{equation*}%
Now, we have
\begin{eqnarray*}
<k>^{s_{1}}\Vert \Box _{k}e^{-t(-\Delta )^{\frac{\alpha }{2}}}f\Vert
_{L^{p}} &\preceq &<k>^{s_{1}-s_{2}}e^{-\frac{t}{2}|k|^{\alpha
}}<k>^{s_{2}}\Vert \Box _{k}f\Vert _{L^{p}} \\
&\preceq &t^{-\frac{1}{\alpha }(s_{1}-s_{2})}<k>^{s_{2}}\Vert \Box
_{k}f\Vert _{L^{p}}.
\end{eqnarray*}%
Taking $l^{q}$ norm in both sides , we obtain (3.14) from the definition of
the modulation space.

Next, we estimate the case $1\leq q_{1}<q_{2}$ and $s_{1}\geq s_{2}$. For
any $\varepsilon >0$, by (2.2) and (3.14), we have
\begin{eqnarray*}
\Vert e^{-t(-\Delta )^{\frac{\alpha }{2}}}f\Vert _{M_{p,q_{1}}^{s_{1}}}
&\preceq &\Vert e^{-t(-\Delta )^{\frac{\alpha }{2}}}f\Vert
_{M_{p,q_{2}}^{s_{1}+\frac{n}{q_{1}}-\frac{n}{q_{2}}+\varepsilon }} \\
&\preceq &(1+t^{-\frac{1}{\alpha }(s_{1}+\frac{n}{q_{1}}-s_{2}-\frac{n}{q_{2}%
}+\varepsilon )})\Vert f\Vert _{M_{p,q_{2}}^{s_{2}}} \\
&=&(1+t^{-\frac{1}{\alpha -\varepsilon _{1}}(\sigma (s_{1},q_{1})-\sigma
(s_{2},q_{2}))})\Vert f\Vert _{M_{p,q_{2}}^{s_{2}}},
\end{eqnarray*}%
where \ \ $\varepsilon _{1}\rightarrow 0+$ \ as \ $\varepsilon \rightarrow
0+.$ \ Notice that the behavior of $1+t^{-\frac{1}{\alpha -\varepsilon }%
[\sigma (s_{1},q_{1})-\sigma (s_{2},q_{2})]}$ likes $t^{-\frac{1}{\alpha
-\varepsilon }[\sigma (s_{1},q_{1})-\sigma (s_{2},q_{2})]}$ when $t$ is
finite. So, by Theorem 1, we can obtain that equation (1.5) is locally
well-posed in $C(0,T;M_{p,q}^{s})$, when $\sigma (s,q)>-\frac{1}{(\alpha
-\varepsilon _{1})(k-1)}$. Since $\ \varepsilon _{1}>0$ \ is arbitrary, we
obtain the conclusion.

\textbf{Proof of Theorem 3 } By the Bejenaru and Tao's conclusion (see Theorem 4 of
\cite{IT}), it suffices to show that the map from $M_{2,q}^{s}$ to $L^{\infty
}([0,T];M_{2,q}^{s})$ defined by
\begin{equation*}
u_{0}\rightarrow \int_{0}^{t}e^{-(t-\tau )\left( -\Delta \right) ^{\frac{%
\alpha }{2}}}(e^{-\tau (-\Delta )^{\frac{\alpha }{2}}}u_{0})^{k}d\tau \eqno%
(3.15)
\end{equation*}%
is discontinuous for $s<-\frac{\alpha }{k-1}$ or $\sigma (s,q)<-\frac{\alpha
}{k-1}$. Actually, if the map is continuous, we will have
\begin{equation*}
\sup\limits_{t\in (0,T)}\Vert \int_{0}^{t}e^{-(t-\tau )\left( -\Delta
\right) ^{\frac{\alpha }{2}}}(e^{\tau \left( -\Delta \right) ^{\frac{\alpha
}{2}}}u_{0})^{k}d\tau \Vert _{M_{2,q}^{s}}\preceq \Vert u_{0}\Vert
_{M_{2,q}^{s}}^{k}.\eqno(3.16)
\end{equation*}%
So, we only need to find a $u_{0}$ such that (3.16) fails.

We first consider the case $s<-\frac{\alpha }{k-1}$. In this case, choose $%
u_{0}$ such that
\begin{equation*}
\mathcal{F}u_{0}=\chi _{N}=\chi (\xi -N\mathbf{e})+\chi (\xi +N\mathbf{e}%
)=\chi _{+}(\xi )+\chi _{-}(\xi )
\end{equation*}%
where $N$ is a large natural number, $\mathbf{e}=(1,1,...,1)$, and $\chi $
is the characteristic function of the cube
\begin{equation*}
E=[-1,1]^{n}.
\end{equation*}%
This $\mathcal{F}u_{0}$ is a non-negative even function. By the choice of \ $%
\mathcal{F}u_{0}$ \ and the definition of the modulation space, using the
Plancerel formula we have
\begin{equation*}
\Vert u_{0}\Vert _{M_{2,q}^{s}}^{k}\preceq N^{ks}.\eqno(3.17)
\end{equation*}%
Now, we estimate
\begin{equation*}
\Vert \int_{0}^{t}e^{-(t-\tau )\left( -\Delta \right) ^{\frac{\alpha }{2}%
}}(e^{-\tau \left( -\Delta \right) ^{\frac{\alpha }{2}}}u_{0})^{k}d\tau
\Vert _{_{M_{2,q}^{s}}}.
\end{equation*}%
By taking $t=\frac{1}{N^{\alpha }},$ we get
\begin{eqnarray*}
&&\Vert \int_{0}^{\frac{1}{N^{\alpha }}}e^{-(\frac{1}{N^{\alpha }}-\tau
)\left( -\Delta \right) ^{\frac{\alpha }{2}}}(e^{-\tau \left( -\Delta
\right) ^{\frac{\alpha }{2}}}u_{0})^{k}d\tau \Vert _{_{M_{2,q}^{s}}}^{q} \\
&~&~~~~~~~~=\sum_{j\in
\mathbb{Z}
^{n}}<j>^{sq}\Vert \Box _{j}\int_{0}^{\frac{1}{N^{\alpha }}}e^{-(\frac{1}{%
N^{\alpha }}-\tau )\left( -\Delta \right) ^{\frac{\alpha }{2}}}(e^{-\tau
\left( -\Delta \right) ^{\frac{\alpha }{2}}}u_{0})^{k}d\tau \Vert
_{L^{2}}^{q}~~~~~~~(3.18)
\end{eqnarray*}%
We denote the convolution of $k$ functions of $\chi _{N}$ by $\chi _{N}\ast
\cdot \cdot \cdot \ast \chi _{N}$. It is easy to find that the cube \ $%
E_{N}=[kN-k,kN+k]^{n}$ \ is a subset of the support of $\chi _{N}\ast \cdot
\cdot \cdot \ast \chi _{N}$. \ Also, notice that
\begin{equation*}
e^{-(\frac{1}{N^{\alpha }}-\tau )|\xi |^{\alpha }}\geq C>0
\end{equation*}%
for $\tau \in \lbrack 0,\frac{1}{N^{\alpha }}]$ and $\xi \in E_{N},$ \ and
that
\begin{equation*}
e^{-\tau |\xi |^{\alpha }}\geq C>0
\end{equation*}%
for $\tau \in \lbrack 0,\frac{1}{N^{\alpha }}]$ and $\xi \in supp\chi _{N}$.
By the Plancerel theorem, we have that, for \ $j=kN\mathbf{e},$%
\begin{eqnarray*}
&&\Vert \Box _{j}\int_{0}^{\frac{1}{N^{\alpha }}}e^{-(\frac{1}{N^{\alpha }}%
-\tau )\left( -\Delta \right) ^{\frac{\alpha }{2}}}(e^{-\tau \left( -\Delta
\right) ^{\frac{\alpha }{2}}}u_{0})^{k}d\tau \Vert _{L^{2}}^{q} \\
&=&\Vert \varphi _{j}(\xi )\int_{0}^{\frac{1}{N^{\alpha }}}e^{-(\frac{1}{%
N^{\alpha }}-\tau )|\xi |^{\frac{\alpha }{2}}}\left\{ (e^{-\tau |\cdot
|^{\alpha }}\chi _{+})\ast \cdot \cdot \cdot \cdot \ast (e^{-\tau |\cdot
|^{\alpha }}\chi _{+})\right\} (\xi )d\tau \Vert _{L^{2}(d\xi )}^{q} \\
&\succeq &C\int_{0}^{\frac{1}{N^{\alpha }}}\Vert (\chi _{+})\ast \cdot \cdot
\cdot \cdot \ast (\chi _{+})\Vert _{L^{2}(E_{N}\bigcap supp\varphi
_{j})}^{q}d\tau .
\end{eqnarray*}
Moreover, because the Lebesgue measure of $E_{N}$ is a constant, we have
that for \ $j=kN\mathbf{e}$
\begin{equation*}
\langle j\rangle^{sq}\Vert \Box _{j}\int_{0}^{\frac{1}{N^{\alpha }}}e^{-(\frac{1}{%
N^{\alpha }}-\tau )\left( -\Delta \right) ^{\frac{\alpha }{2}}}(e^{-\tau
\left( -\Delta \right) ^{\frac{\alpha }{2}}}u_{0})^{k}d\tau \Vert
_{L^{2}}^{q}\succeq N^{\left( s-\alpha \right) q}.
\end{equation*}%
It leads to the inequality
\begin{equation*}
\Vert \int_{0}^{\frac{1}{N^{\alpha }}}e^{-(\frac{1}{N^{\alpha }}-\tau
)\left( -\Delta \right) ^{\frac{\alpha }{2}}}(e^{-\tau \left( -\Delta
\right) ^{\frac{\alpha }{2}}}u_{0})^{k}d\tau \Vert _{_{M_{2,q}^{s}}}\geq
CN^{s-\alpha }\eqno(3.19)
\end{equation*}%
which contradicts to (3.16) and (3.17). So equation (1.5) is ill-posed in $%
M_{2,q}^{s}$ when $s<-\frac{\alpha }{k-1}.$

Now, we consider the case $\sigma (s,q)<-\frac{\alpha }{k-1}$. For
convenience, we let
\begin{equation*}
\mathcal{F}u_{0}=\chi _{N}^{\ast }=\chi (\frac{1}{N}(\xi -100kN\mathbf{e)}),
\end{equation*}%
where \ $\chi $ \ is the characteristic function of the set \ $E=[-1,1]^{n}.$
If we want $u_{0}$ to be a real function, we can make an even extension just
like what we did in the previous case, the result should be the same. So, by
the Plancherel theorem and the definition of the modulation spaces, we have
\begin{equation*}
\Vert u_{0}\Vert _{M_{2,_{q}}^{s}}^{k}\preceq N^{k(s+\frac{n}{q})}\eqno%
(3.20).
\end{equation*}%
On the other hand, choosing $t=N^{-\alpha }$ again, by the similar method as
we did previously, when \ $j$ \ is the center (or very close to the center
of the support of \ $\chi _{N}^{\ast }\ast \cdot \cdot \cdot \cdot \ast \chi
_{N}^{\ast },$ \ \ we have
\begin{eqnarray*}
&&\Vert \Box _{j}\int_{0}^{\frac{1}{N^{\alpha }}}e^{-(\frac{1}{N^{\alpha }}%
-\tau )\left( -\Delta \right) ^{\frac{\alpha }{2}}}(e^{-\tau \left( -\Delta
\right) ^{\frac{\alpha }{2}}}u_{0})^{k}d\tau \Vert _{L^{2}}^{q} \\
&\geq &C\Vert \int_{0}^{\frac{1}{N^{\alpha }}}e^{-(\frac{1}{N^{\alpha }}%
-\tau )|\xi |^{\frac{\alpha }{2}}}(e^{-\tau |\cdot |^{\alpha }}\chi
_{N}^{\ast })\ast \cdot \cdot \cdot \cdot \ast (e^{-\tau |\cdot |^{\alpha
}}\chi _{N}^{\ast })d\tau \Vert _{L^{2}(supp\varphi _{j})}^{q} \\
&\geq &C\Vert \int_{0}^{\frac{1}{N^{\alpha }}}(\chi _{N}^{\ast })\ast \cdot
\cdot \cdot \cdot \ast (\chi _{N}^{\ast })d\tau \Vert _{L^{2}(supp\varphi
_{j})}^{q} \\
&\geq &CN^{-\alpha }\Vert (\chi _{N}^{\ast })\ast \cdot \cdot \cdot \cdot
\ast (\chi _{N}^{\ast })\Vert _{L^{2}(supp\varphi _{j})}^{q}.
\end{eqnarray*}%
Notice that the Lebesgue measure of supp $\chi _{N}^{\ast }$ is $N^{n}$
times that of supp $\Box _{i}$, and $(\chi _{N}^{\ast }\ast \cdot \cdot
\cdot \cdot \ast \chi _{N}^{\ast })$ is constructed by $(k-1)$ convolutions.
Therefore, we have
\begin{equation*}
\Vert \chi _{N}^{\ast }\ast \cdot \cdot \cdot \cdot \ast \chi _{N}^{\ast
}\Vert _{L^{2}(supp\varphi _{i})}\succeq N^{(k-1)n}.\eqno(3.21)
\end{equation*}%
Moreover, the support of $\chi _{N}^{\ast }\ast \cdot \cdot \cdot \cdot \ast
\chi _{N}^{\ast }$ \ is the cube
\begin{equation*}
E_{N}^{\ast }:=[100k^{2}N-kN,100^{2}kN+kN]^{n}
\end{equation*}%
of Lebesgue measure $\left( 2kN\right) ^{n}$. Therefore, the number of
summands, in right side of (3.18) is $CN^{n}$ for some constant \ $C>0.$ We
now obtain
\begin{equation*}
\Vert \int_{0}^{\frac{1}{N^{\alpha }}}e^{(\frac{1}{N^{\alpha }}-\tau )\Delta
^{\frac{\alpha }{2}}}(e^{\tau \Delta ^{\frac{\alpha }{2}}}u_{0})^{k}d\tau
\Vert _{_{M_{2,q}^{s}}}\geq CN^{s+\frac{n}{q}+(k-1)n-\alpha }.\eqno(3.22)
\end{equation*}%
The last inequality gives a contradiction to (3.16) and (3.20). So equation
(1.5) is ill-posed in $M_{2,q}^{s}$ when $\sigma (s,q)<-\frac{\alpha }{k-1}$%
. This completes the proof of Theorem 3.

\bigskip

The proofs for Corollary 1 and Corollary 2 are similar to the above proof.
We leave them to the reader.

\end{document}